\documentclass{cmslatex}
\usepackage{latexsym,amssymb,enumerate,amsmath,graphicx}

\renewcommand{\theequation}{\arabic{section}.\arabic{equation}}

\newtheorem{thm}{Theorem}[section]
\newtheorem{lem}[thm]{Lemma}
\newtheorem{cor}[thm]{Corollary}

\newcommand{\R}{{\mathbb R}}
\newcommand{\D}{{\mathrm D}}
\newcommand{\DD}{{\mathrm D}^2}
\newcommand{\Dg}{\Delta_g}
\newcommand{\ig}[1]{\int_\Omega{#1}\;d\gamma}
\newcommand{\dig}[1]{\int_{\partial\Omega}{#1}\;g\,d{\mathcal H}^{d-1}}
\newcommand{\Ep}{{\mathcal E}_p}
\newcommand{\Ip}{{\mathcal I}_p}
\newcommand{\Kp}{{\mathcal K}_p}
\newcommand{\Epm}{{\mathcal E}_{m,p}}
\newcommand{\Ipm}{{\mathcal I}_{m,p}}
\newcommand{\Kpm}{{\mathcal K}_{m,p}}
\newcommand{\be}[1]{\begin{equation}\label{#1}}
\newcommand{\ee}{\end{equation}}
\newcommand{\eqn}[1]{(\ref {#1})}
\newcommand{\prf}{\par\noindent{\sl Proof.\/} }
\newcommand{\finprf}{\unskip\null\hfill$\square$\vskip 0.3cm}

\renewcommand{\(}{\left(}
\renewcommand{\)}{\right)}
\newcommand{\nrm}[2]{\|#1\|_{L^{#2}(\Omega,d\gamma)}}
\newcommand{\email}[1]{{\small E-mail: {\textsf {#1}}}}
\newcommand{\urladdr}[1]{{\small URL: {\textsf {#1}}}}
\newcommand{\address}[1]{\thanks{#1}}

\begin{document}
\title
{On the Bakry-Emery criterion for linear diffusions and weighted porous media equations}
\author{J. Dolbeault, B. Nazaret, and G. Savar\'e}
\date{\today}
\maketitle
\thispagestyle{empty}

\begin{abstract} The goal of this paper is to give a non-local sufficient condition for generalized Poincar\'e inequalities, which extends the well-known Bakry-Emery condition. Such generalized Poincar\'e inequalities have been introduced by W. Beckner in the gaussian case and provide, along the Ornstein-Uhlenbeck flow, the exponential decay of some generalized entropies which interpolate between the $L^2$ norm and the usual entropy. Our criterion improves on results which, for instance, can be deduced from the Bakry-Emery criterion and Holley-Stroock type perturbation results. In a second step, we apply the same strategy to non-linear equations of porous media type. This provides new interpolation inequalities and decay estimates for the solutions of the evolution problem. The criterion is again a non-local condition based on the positivity of the lowest eigenvalue of a Schr\"odinger operator. In both cases, we relate the Fisher information with its time derivative. Since the resulting criterion is non-local, it is better adapted to potentials with, for instance, a non-quadratic growth at infinity, or to unbounded perturbations of the potential. \end{abstract}

\bigskip

The Bakry-Emery method \cite{Bakry-Emery85} has been extremely successful to establish logarithmic Sobolev inequalities \cite{Gross75}. It is also known to apply very well to the proof of Poincar\'e inequalities, and inequalities which interpolate between Poincar\'e and logarithmic Sobolev inequalities \cite{MR0954373,Arnold-Markowich-Toscani-Unterreiter01,MR1845806,MR2081075}, the so-called ``generalized Poincar\'e inequalities.'' The first paper on such inequalities has been written by W. Beckner in \cite{MR0954373} in the case of a gaussian measure. In \cite{Arnold-Markowich-Toscani-Unterreiter01}, the inequalities have been slightly generalized by taking into account general $x$ dependent diffusions and by considering ``convex entropies'' based on a convex function $\psi$ satisfying the additional admissibility condition: $\psi''\,\psi^{(iv)}\geq 2\,(\psi''')^2$. The proof is based on the entropy -- entropy production method and the Bakry-Emery condition, \cite{Bakry-Emery85}. The result has been slightly improved in~\cite{MR2152502}, thus providing ``refined inequalities'' in the case $\psi(s)=(s^p-1-p\,(s-1))/(p-1)$, $p\in (1,2)$. Other considerations on ``generalized Poincar\'e inequalities'' can be found in \cite{MR1845806,MR1796718,MR2081075,ABD}.

The Bakry-Emery condition is a sufficient, local condition, which relies on a uniform strict log-concavity of the measure. Using perturbation techniques, it is possible to relax such a strict assumption to some extent, see \cite{MR893137,MR2081075}. When Poincar\'e and logarithmic Sobolev inequalities are known to hold simultaneously, further interpolation inequalities can also be established: see for instance \cite{MR1796718,MR1038450,Bartier-Dolbeault,ABD}.

In case of generalized Poincar\'e inequalities, many techniques which are available for Poincar\'e inequalities and spectral gap approaches can be adapted, and more flexibility is available than for logarithmic Sobolev inequalities. See \cite{Arnold-Markowich-Toscani-Unterreiter01,MR1796718,MR2152502,MR2081075}. Our first main result is that {\em the Bakry-Emery condition can be replaced by a non-local condition which amounts to assume the positivity of the first eigenvalue of a Schr\"odinger operator}. This criterion gives only a sufficient condition, but improves on a very popular method, as will be illustrated by two examples.

\medskip The second main result of this paper is concerned with interpolation inequalities and decay estimates for the solutions of a weighted porous media equation. In \cite{DelPino-Dolbeault01a}, del Pino and the first author proved by variational techniques that sharp constants in some Gagliardo-Nirenberg interpolation were equivalent to sharp exponential decay rates of the generalized entropy of Ralston and Newman \cite{MR760592,MR760591} along the flow associated to the porous media equation in the euclidean space. A similar result was achieved in \cite{MR2001j:35155} by an entropy -- entropy production method, and later extended in \cite{Carrillo-Juengel-Markowich-Toscani-Unterreiter01,MR2003b:35118,Carrillo-Vazquez03}. Also see \cite{MR2213477}, and \cite{MR2065020} for a review. A very nice interpretation in terms of a gradient flow was done first in the limit case of the heat equation in \cite{MR1617171}, and then in the porous media case in \cite{Otto01}. Also see \cite{MR2129498,MR1982656,MR2211152,MR2192294,DNS} for further results based on Wasserstein's distance and mass transportation theory. Concerning interpolation inequalities, weights and asymptotic behaviour, we also have to quote \cite{BBDGV-Cras,Daskalopoulos-Sesum2006,BBDGV2} for some recent results.

In connection with probability theory, entropy -- entropy production methods have been successfully applied in a quite general framework of Riemanian manifolds: see \cite{Bakry02,MR2178944}. However, our approach is more related to a series of attempts which have been made to establish rates of decay for solutions of diffusion equations of second and higher order: see \cite{MR2219352,MR2224869,MR2209292}. More specifically, in \cite{MR2224869}, the question of how generalized Poincar\'e inequalities can be adapted to solutions of the porous media equation on the one-dimensional periodic interval was studied by direct methods, and later extended in \cite{0661} using capacity -- measure tools. The main originality of our second result is probably that {\em it relates interpolation inequalities to a linear eigenvalue problem, and therefore provides a non-local criterion}, which is again an extension of the ideas of Bakry-Emery. This criterion is the same as in the linear case.

We will not discuss here the question of the regularity of the solutions and the well-posedness of the Cauchy problem. The interested reader will find in \cite{MR2286292} an up-to-date overview of results concerning the porous media equation, and further references therein. In this paper we will only focuse on the derivation of functional inequalities in connection with decay estimates for entropy and entropy production estimates of smooth positive solutions of diffusion equations. We shall assume that the solutions are smooth and behave adequately under integrations by parts. The functional inequalities then hold in the functional spaces for which all integral quantities are well defined, as can be shown by a standard density argument. Decay estimates for the solutions of the evolution problems should rather be considered as an {\sl a priori\/} estimate which defines a special class of solutions to the Cauchy problem. Such estimates are again standard in the framework of an approximation scheme.

\section{A nonlocal criterion for the Ornstein-Uhlenbeck flow}

\subsection{Notations and preliminary computations}

Consider on a domain $\Omega\subset\R^d$ the measure $d\gamma=g\,dx$ for some measurable function $g=e^{-F}$ and define the Ornstein-Uhlenbeck operator $\Dg$ by $\Dg v:=\Delta v-\D F\cdot\D v$. The fundamental property of this operator is that, for any $v\in H^1_0(\Omega,d\gamma)$,
\[
\ig{|\D v|^2}=-\ig{v\,\Dg v}\;.
\]
For simplicity, we assume that $\partial\Omega$ is of class $C^1$ and that $g$ is smooth enough so that the trace of $g$ is well defined on $\partial\Omega$. We shall therefore consider the measure $g\,d{\mathcal H}^{d-1}$ induced by $d\gamma$ on $\partial\Omega$ and assume that regularity is such that integrations by parts can be performed without further precautions. Here $H^1_0(\Omega,d\gamma)$ is the space $\{w\in W_{\rm loc}^{1,1}(\R^d,dx)\cap L^2(\Omega,d\gamma)\,:\,\|u\|_{L^2(\partial\Omega,d\gamma)}=0\;\mbox{and}\;\|Du\|_{L^2(\Omega,d\gamma)}<\infty\}$. Let $v$ be a positive smooth solution of the equation
\be{Eqn:Heat}
v_t=\Dg v\quad x\in\Omega\,,\; t\in\R^+\,,
\ee
satisfying homogeneous Neumann boundary conditions
\be{Eqn:bc}
\nabla v\cdot n=0\quad x\in\partial\Omega\,,\; t\in\R^+\,,
\ee
where $n=n(x)$ denotes the unit outgoing normal vector at $x\in\partial\Omega$. Define for $p\in (1,2]$ the functionals
\begin{eqnarray*}
\Ep(t)&:=&\frac 1{p-1}\ig{\Big[v^p-1-p\,(v-1)\Big]}\;,\\
\Ip(t)&:=&\frac 4p\ig{|\D s|^2}\;,\\
\Kp(t)&:=&\ig{|\Dg s|^2}+\alpha\!\ig{\Dg s\,\frac{|\D s|^2}s}\;,
\end{eqnarray*}
with $s:=v^{p/2}$ and $\alpha:=(2-p)/p$. The limit case $p=1$ is also covered with $v=s^2$
\[
\mathcal E_1(t):=\ig{\Big[v\,\log v-(v-1)\Big]}\;,
\]
and $\mathcal I_1$ and $\mathcal K_1$ as above. We shall say that $\Ep$, $p\in [1,2]$, is the {\sl generalized entropy,\/} or simply the {\sl entropy,\/} since the case $p>1$ is clearly a generalization of the case $p=1$. Written in terms of $s$, for $p>1$ the entropy is
\[
\Ep=\frac 1{p-1}\ig{\Big[s^2-1-p\,(s^{2/p}-1)\Big]}
\]
and the evolution is governed by
\be{Eqn:Heat2}
s_t=\Dg s+\alpha\,\frac{|\D s|^2}s\;.
\ee
Using \eqn{Eqn:Heat} and \eqn{Eqn:Heat2}, a computation shows that
\begin{eqnarray*}
\frac d{dt}\Ep(t)&:=&-\Ip(t)\;,\\
\frac d{dt}\Ip(t)&:=&-\frac 8p\,\Kp(t)\;.
\end{eqnarray*}

\bigskip The next step consists in establishing a functional inequality relating $\Ip$ and $\Kp$. Using the commutation relation
\[
[\D,\Dg]\,s=-\DD F\,\D s
\]
and integration by parts, one gets
\[
\ig{|\Dg s|^2}=\ig{|\DD s|^2}+\ig{\DD F\,\D s\cdot\D s}-\sum_{i,j=1}^d\dig{\partial_{ij}^2s\,\partial_i s\, n_j}
\]
if Neumann boundary conditions hold.
\medskip

\begin{lem}\label{Lem:Convex} With the above notations, if $\Omega$ is convex, then
\[
-\sum_{i,j=1}^d\dig{\partial_{ij}^2s\,\partial_i s\, n_j}\geq 0\;.
\]
\end{lem}

\medskip
For a proof of this result, see \cite[Lemma 5.2]{GTS} or \cite{MR775683}.

Hence, if $\Omega$ is convex, we have
\be{Ineq:IPP1}
\ig{|\Dg s|^2}\geq\ig{|\DD s|^2}+\ig{\DD F\,\D s\cdot\D s}
\ee

With the notation $z:=\sqrt s$, we get
\[
\D z=\frac{\D s}{2\sqrt s}
\]
and, using integration by parts again, since $2\,\DD s:\D z\otimes\D z=\D\(|\D z|^2\):\D z\,$,
\be{Ineq:IPP2}
4\ig{|D z|^4}=\ig{\Dg s\,|\D z|^2}+2\ig{\DD s\,:\,\D z\otimes\D z}\;.
\ee
Hence by \eqn{Ineq:IPP1} and \eqn{Ineq:IPP2} we have
\begin{eqnarray*}
\Kp&=&\ig{|\Dg s|^2}+4\,\alpha\!\ig{\Dg s\,|\D z|^2}\\
&\geq&\ig{|\DD s|^2}+\ig{\DD F\,\D s\cdot\D s}\\
&&\qquad +\; 4^2\,\alpha\!\ig{|D z|^4}-2\cdot 4\,\alpha\!\ig{\DD s\,:\,\D z\otimes\D z}\;.
\end{eqnarray*}
Expanding the square, we get
\[
\ig{\left|\DD s-4\,\D z\otimes\D z\right|^2}=\ig{|\DD s|^2}-2\cdot 4\ig{\DD s\,:\,\D z\otimes\D z}+4^2\ig{|D z|^4}\;,
\]
which gives the identity
\[
4^2\,\alpha\!\ig{|D z|^4}-2\cdot 4\,\alpha\!\!\ig{\DD s\,:\,\D z\otimes\D z}=\alpha\!\!\ig{\left|\DD s-4\,\D z\otimes\D z\right|^2}-\alpha\!\!\!\ig{|\DD s|^2}
\]
and hence
\be{Ineq:Decay}
\Kp\geq (1-\alpha)\ig{|\DD s|^2}+\ig{\DD F\,\D s\cdot\D s}\;.
\ee
For $p=1$, only the last term of the right hand side is left. For $p>1$, the coefficient $1-\alpha=2\,(p-1)/p$ is also positive, which allows to state a nonlocal result.

\subsection{A nonlocal result}\label{Sec:NonLoc}
Let
\[
V(x):=\inf_{\xi\in S^{d-1}}\big(\DD F(x)\,\xi,\xi\big)
\]
and define
\[
\lambda_1(p):=\inf_{w\in H^1(\Omega,d\gamma)\setminus\{0\}}\frac{\ig{\left(2\,\frac{p-1}p\,|\D w|^2+V\,|w|^2 \right)}}{\ig{|w|^2}}\;.
\]

Notice that the above criterion is given with respect to the measure $d\gamma$ and $H^1(\Omega,d\gamma)$ is defined as the space $\{w\in W_{\rm loc}^{1,1}(\R^d,dx)\cap L^2(\Omega,d\gamma)\,:\,\|Du\|_{L^2(\Omega,d\gamma)}<\infty\}$. If we assume that $F\in C^2(\bar\Omega)$ is such that
\[
DF\cdot n\geq 0\quad\mbox{on}\quad\partial\Omega
\]
(or $\Omega=\R^d$), and consider
\[
u:=w\,e^{-F/2}\,,
\]
a simple integration by parts provides a sufficient condition for the positivity of $\lambda_1(p)$ in terms of the positivity of the ground state level of a Schr\"odinger equation in the usual sense:
\[
\lambda_1(p) \geq 2\,\frac{p-1}p\,\inf_{w\in H^1(\Omega,dx)\setminus\{0\}}\frac{\int_\Omega{\left[|\D u|^2+\(\frac p{2(p-1)}\,V+\frac 14\,|\D F|^2-\frac 12\,\Delta F\)\,|u|^2 \right]}\,dx}{\int_\Omega{|u|^2}\,dx}\;.
\]
Notice that the above condition is only a sufficient condition for relating $\Ip$ and $\Kp$ through \eqn{Ineq:Decay}. We are indeed reducing the problem to a scalar eigenvalue problem and dropping the fact that $\D s\cdot n=0$ on $\partial\Omega$.
\medskip

\begin{thm}\label{Thm:Global} Assume that $F\in C^2(\bar\Omega)$, $e^{-F}\in L^1(\Omega)$, where $\Omega$ is a convex domain in~$\R^d$. Consider a smooth solution of \eqref{Eqn:Heat} with Neumann boundary conditions on~$\partial\Omega$ and an initial datum $v_0\in L^1_+(\Omega)$ such that $\Ep(0)$ is finite for some $p\in [1,2]$. With the above notations, if
\be{Hyp:OU}
\lambda_1(p)>0\;,
\ee
then the solution of \eqn{Eqn:Heat} is such that, for any $t\geq 0$,
\begin{eqnarray*}
&&\Ip(t)\leq \Ip(0)\,e^{-2\,\lambda_1(p)\,t},\\
&&\Ep(t)\leq \Ep(0)\,e^{-2\,\lambda_1(p)\,t}.
\end{eqnarray*}\end{thm}

\medskip
\prf The first inequality is a consequence of \eqn{Ineq:Decay}. We may indeed write
\begin{eqnarray*}
-\frac p8\,\frac d{dt}\Ip(t)=
\Kp&\geq& (1-\alpha)\ig{|\DD s|^2}+\ig{V(x)\,|\D s|^2}\\
&\geq& \lambda_1(p)\ig{|\D s|^2}=2\,\lambda_1(p)\,\frac p8\,\Ip(t)
\end{eqnarray*}
and use a Gronwall lemma to prove the first inequality. The second inequality then follows from
\[
\frac d{dt}\Big[\Ip(t)-2\,\lambda_1(p)\,\Ep(t)\Big]\leq 0\;,
\]
the fact that
\[
\lim_{t\to\infty}\Ip(t)=0
\]
and, as a consequence,
\[
\lim_{t\to\infty}\nrm{v(t)-\bar v}1=0\;,\quad\mbox{where}\quad \bar v=\frac{\ig{v_0}}{\ig{}}\;,
\]
so that
\be{Eqn:Zero}
\lim_{t\to\infty}\Ep(t)=0\;.
\ee

An integration on $(t,\infty)$ then gives
\be{Ineq:Gronwall}
\frac d{dt}\Ep(t)-2\,\lambda_1(p)\,\Ep(t)=\Ip(t)-2\,\lambda_1(p)\,\Ep(t)\geq \lim_{t\to\infty}\(\Ip(t)-2\,\lambda_1(p)\,\Ep(t)\)=0\;,
\ee
and one again concludes with a Gronwall lemma.\finprf

As a standard consequence, this proves a generalized Poincar\'e inequality. Notice that the result was already known at least in case $p=2$, see, {\sl e.g.,\/} \cite[Proposition 5.5.4]{MR1845806}. How the so-called {\sl super-integral\/} criterion for the logarithmic Sobolev inequality, see \cite[Proposition 5.5.6]{MR1845806}, is related to our approach in the limit $p\to 1$ is still an open question. The boundedness of $d\gamma$, that is $e^{-F}\in L^1(\Omega,dx)$, is used only to obtain~\eqn{Eqn:Zero}. Whenever \eqn{Eqn:Zero} can be established by other means, the method also adapts to unbounded measures.

\medskip

\begin{cor}\label{Cor:Beckner} Assume that $F\in C^2(\bar\Omega)$, $e^{-F}\in L^1(\Omega)$, where $\Omega$ is a convex domain in $\R^d$. If $\lambda_1(p)$ is positive, then
\be{Eqn:GPI}
\frac 1{p-1}\left[\;\ig{|u|^2}-\(\ig{|u|^\frac 2p}\)^p\;\right]\leq \frac 2{\lambda_1(p)}\ig{|\nabla u|^2}\quad\forall\; u\in H^1(\Omega,d\gamma)\;.
\ee\end{cor}

\medskip
\begin{proof} As a straightforward consequence of the proof of Theorem~\ref{Thm:Global}, namely according to~\eqref{Ineq:Gronwall} taken at $t=0$, we have
\[
\Ep(t=0)\leq \frac 1{2\,\lambda_1(p)}\,\Ip(t=0)\;,
\]
which means
\[
\frac 1{p-1}\ig{\left[v_0^p-1-p\,(v_0-1)\right]}\leq \frac 1{2\,\lambda_1(p)}\,\frac 4p\ig{|\D v_0^{p/2}|^2}
\]
for any initial value $v_0$ of the Cauchy problem associated to \eqref{Eqn:Heat}. Let $u:=\mu\,v_0^{p/2}$. The inequality now takes the form
\[
\frac 1{p-1}\ig{\left[u^2-\mu^2-p\,\(u^\frac 2p\,\mu^{2\frac{p-1}p}-\mu^2\)\right]}\leq \frac 1{2\,\lambda_1(p)}\,\frac 4p\ig{|\D u|^2}\;.
\]
If we observe that the left hand side, as a function of $\mu$, can be written
\[
h(\mu)+\frac 1{p-1}\ig{|u|^2}
\]
where $h(\mu):=\mu^2-\frac p{p-1}\,|\bar u|^\frac 2p\,\mu^{2\frac{p-1}p}$ and $\bar u=\(\ig{u^\frac 2p}\)^\frac p2$, an optimization of $h(\mu)$ in terms of $\mu$ shows that
\[
h(\mu)\geq h(\bar u)=-\frac 1{p-1}\,|\bar u|^2\quad\forall\;\mu>0\;,
\]
which completes the proof. \end{proof}

One can notice that under the {\sl local\/} (uniform) strict convexity condition
\[
\inf_{\xi\in S^{d-1}}\big(\DD F(x)\,\xi,\xi\big)=:\lambda_1^\infty>0\;,
\]
one recovers the usual results corresponding to the Bakry-Emery criterion~\cite{Bakry-Emery85}, since in such a case $\lambda_1(p)\geq \lambda_1^\infty$ for any $p\in [1,2]$. Our approach does not cover all cases. If~$\Omega$ is bounded and $F$ is constant, the infimum on $H^1(\Omega,d\gamma)\setminus\{0\}$ in the definition of $\lambda_1(p)$ is zero, which is achieved by constant functions. In such a case, it is therefore crucial to take into account the condition $\D s\cdot n=0$ on $\partial\Omega$.

If $V$ is nonnegative, $\lambda_1(p)\geq (p-1)\,\lambda_1(2)$. The inequality
\be{Ineq:BecknerPoincare}
\frac 1{p-1}\left[\;\ig{|u|^2}-\(\ig{|u|^\frac 2p}\)^p\;\right]\leq \frac 2{(p-1)\,\lambda_1(2)}\ig{|\nabla u|^2}\quad\forall\; u\in H^1(\Omega,d\gamma)\;.
\ee
also holds for any $p\in(1,2)$ if it holds for $p=2$, even if $V$ changes sign. Using H\"older's inequality and the fact that $d\gamma$ is a probability measure, it is indeed not difficult to see that
\[
\ig{|u|^2}-\(\ig{|u|^\frac 2p}\)^p\leq\ig{|u|^2}-\(\ig{|u|}\)^2\,.
\]
All results which can be achieved for $p=2$ can therefore be extended to $p\in(1,2)$, with a constant that can be bounded by the same constant as for $p=2$, up to a factor $1/(p-1)$. However, requiring that $\lambda_1(p)>0$ gives much better estimates: first of all, it allows to pass to the limit case $p=1$ in the standard case of the Bakry-Emery method, that is when $\lambda_1^\infty>0$, with a better constant for any $p>1$. Moreover, it allows do directly deal with unbounded perturbations or limit cases, as illustrated by the two following examples.

\subsection{Two examples}

\subsubsection*{1) Unbounded perturbations}

Let $F(x)=\frac 12\,|x|^2+\varepsilon\,\log|x|$, $\Omega=\R^d$ for some $d\geq 3$. The problem cannot be reduced to the harmonic potential case by the Holley-Stroock perturbation lemma, see \cite{MR893137,Arnold-Markowich-Toscani-Unterreiter01,BDK2007}, since at the origin $x=0$ and as $|x|\to\infty$, the perturbation $\varepsilon\,\log|x|=F(x)-\frac 12\,|x|^2$ is clearly unbounded. Assume that $d\geq 3$. The measure $d\gamma$ is bounded for any $\varepsilon\in (0,d)$. With $\nu:=\frac p{2\,(p-1)}\geq 1$, using Hardy's inequality and the expression of the first eigenvalue for the quantum harmonic oscillator, we can write that
\begin{multline}\label{Ineq:Unbounded}
\int_\Omega{\left[\,|\D u|^2+\big({\textstyle{\nu\,V+\frac 14\,|\D F|^2-\frac 12\,\Delta F}}\big)\,|u|^2\,\right]}\,dx\\
=(1-\mathsf a^2)\kern-3pt\underbrace{\int_\Omega{\left[|\D u|^2-{\textstyle{\frac{(d-2)^2}{4\,|x|^2}}}\,|u|^2\right]}\,dx}_{\geq 0}+\kern-3pt\underbrace{\int_\Omega{\left[\mathsf a^2\,|\D u|^2+{\textstyle{\frac{|x|^2}4}}\,|u|^2\right]}\,dx}_{\geq \frac d2\,\mathsf a\int_\Omega |u|^2\,dx}+\(\nu-{\textstyle{\frac{d-\varepsilon}2}}\)\kern-3pt\int_\Omega\!|u|^2\,dx
\end{multline}
with $\mathsf a$ such that $(1-\mathsf a^2)\,\frac{(d-2)^2}4=\frac\varepsilon{4}\(2\,\mathsf b-\varepsilon\)$ and $\mathsf b:=2\,\nu+d-2$. This can be done if and only if $0\leq1-\mathsf a^2=\frac\varepsilon{(d-2)^2}\,\(2\,\mathsf b-\varepsilon\)\leq 1$, that is
\be{CdtZerOne}
\varepsilon\;\leq\;\mathsf b-\sqrt{\mathsf b^2-(d-2)^2}\;.
\ee
The case $\varepsilon\geq\mathsf b+\sqrt{\mathsf b^2-(d-2)^2}$ has indeed to be excluded because $\mathsf b\geq d>\varepsilon$. For $\nu>\frac d2$, {\sl i.e.\/} $p<\frac d{d-1}$, the last term of the right hand side in \eqn{Ineq:Unbounded} is positive, which proves that $\lambda_1(p)>0$ under condition \eqn{CdtZerOne}. As a special case, this means that an upper bound for $\varepsilon$ is of the order of $\frac{(d-2)^2}{4\,\nu}=\frac{(d-2)^2}{2\,p}\,(p-1)$ as $p\to 1$. For $\nu\leq\frac d2$, we can evaluate $\lambda_1(p)$ from below by $\frac d2\,\mathsf a+\nu-\frac{d-\varepsilon}2$, which is positive for $\varepsilon>0$, small enough.

\subsubsection*{2) Semi-classical estimates}

The criterion of Theorem~\ref{Thm:Global} amounts to look for the ground state of
\[
-\Delta u+\(\nu\,V+\frac 14\,|\D F|^2-\frac 12\,\Delta F\)\,u=\nu\,\lambda_1(p)\,u
\]
with $\nu:=\frac p{2\,(p-1)}\geq 1$ as above, and to determine $\lambda_1(p)$. Let $F(x)=\frac 1\beta\,|x|^\beta$, $\beta\in(1,2)$, $\Omega=\R$. In such a case, the problem is reduced to
\[
-u''+\left[\(\textstyle{\nu-\frac 12}\)\,\(\textstyle{\beta-1}\)\,|x|^{\beta-2}+\frac 14\,|x|^{2(\beta-1)}\right]\,u=\nu\,\lambda_1(p)\,u\;.
\]
It readily follows that the criterion of Theorem~\ref{Thm:Global} covers the range $1<\beta<2$, which is not the case of the usual Bakry-Emery criterion, and a semi-classical expansion shows that
\[
\lambda_1(p)\sim\kappa_\beta\,\(\textstyle{\nu-\frac 12}\)^{\frac 2\beta(\beta-1)}\nu^{-1}=O\((p-1)^\frac{2-\beta}\beta\)\quad\mbox{as}\quad p\to 1_+
\]
for some positive constant $\kappa_\beta$ which depends only on $\beta$. Such an estimate is definitely better than the one that can be deduced from \eqn{Ineq:BecknerPoincare}.

\subsection{A refined result in the limit case}\label{Sec:Refined}
Consider the case for which $\lambda_1(p)=0$, but assume that for some $\varepsilon>0$,
\[
\inf_{w\in H^1(\Omega,d\gamma)}\frac{\ig{\left((1-\alpha(1+\varepsilon))\,|\D w|^2+V\,|w|^2 \right)}}{\ig{|w|^2}}=0\;.
\]
Recall that
\begin{eqnarray*}
\Kp&\geq&\ig{|\DD s|^2}+\ig{\DD F\,\D s\cdot\D s}\\
&&\qquad +\; 4^2\,\alpha\!\ig{|D z|^4}-2\cdot 4\,\alpha\!\ig{\DD s\,:\,\D z\otimes\D z}\;.
\end{eqnarray*}
Expanding the square, we get
\begin{eqnarray*}
&&\ig{\left|\sqrt{1+\varepsilon}\;\DD s-\frac 4{\sqrt{1+\varepsilon}}\,\D z\otimes\D z\right|^2}\\
&&=(1+\varepsilon)\ig{|\DD s|^2}-2\cdot 4\ig{\DD s\,:\,\D z\otimes\D z}+\frac{4^2}{1+\varepsilon}\ig{|D z|^4}\;,
\end{eqnarray*}
which gives the identity
\begin{eqnarray*}
&&4^2\,\alpha\!\ig{|D z|^4}-2\cdot 4\,\alpha\!\ig{\DD s\,:\,\D z\otimes\D z}\\
&&=\alpha\!\ig{\left|\sqrt{1+\varepsilon}\,\DD s-\frac 4{\sqrt{1+\varepsilon}}\,\D z\otimes\D z\right|^2}-\alpha\,(1+\varepsilon)\ig{|\DD s|^2}\\
&&\qquad +\;4^2\,\alpha\,\frac \varepsilon{1+\varepsilon}\ig{|D z|^4}\;.
\end{eqnarray*}
Summarizing, we have
\[
\Kp\geq 4^2\,\frac{\alpha\,\varepsilon}{1+\varepsilon}\ig{|D z|^4}\;.
\]
For $p\in(1,2)$, {\sl i.e.\/} $\alpha\in(0,1)$, we follow the method of \cite{MR2152502}. Assume that $d\gamma$ is a probability measure and normalize $v$ in $L^1(\Omega,d\gamma)$ by $\ig v=1$. Then using $\ig{|s|^2}=1+(p-1)\,\Ep$, we can interpolate $\Ip$ with a Cauchy-Schwarz inequality:
\[
\Big(p\,\Ip\Big)^2=4^4\left(\ig{s\,|\D z|^2}\right)^2\leq 4^4\,\Big(1+(p-1)\,\Ep\Big)\ig{|\D z|^4}\;,
\]
which gives
\[
\frac{1+\varepsilon}{\alpha\,\varepsilon}\,\Kp\geq 4^2\ig{|\D z|^4}\geq \Big(\frac p4\Big)^2\,\frac{\,\Ip^2}{1+(p-1)\,\Ep}\;,
\]
which is the functional inequality which replaces \eqn {Eqn:GPI} and gives, a consequence,
\[
-\frac{\frac d{dt}\Ip}{\Ip^2}\geq \frac{\alpha\,\varepsilon}{1+\varepsilon}\,\frac p2\,\frac 1{1+(p-1)\,\Ep(t)}\geq \frac{\alpha\,\varepsilon}{1+\varepsilon}\,\frac p2\,\frac 1{1+(p-1)\,\Ep(0)}=:\kappa\;.
\]
This proves the following result.
\medskip

\begin{thm}\label{Thm:Refined} Assume that $F\in C^2(\bar\Omega)$, $e^{-F}\in L^1(\Omega)$, where $\Omega$ is a convex domain in~$\R^d$, and $\ig{}=1$ so that $d\gamma$ is a probability mesure. Consider a smooth solution of~\eqref{Eqn:Heat} with Neumann boundary conditions on $\partial\Omega$ and an initial datum $v_0\in L^1_+(\Omega)$ such that $\Ep(0)$ is finite for some $p\in [1,2]$. With the above notations, if $\lambda_1(p)=0$, then
\[
\Ip(t)\leq \frac{\Ip(0)}{1+\kappa\,\Ip(0)\,t}\quad\forall\;t\geq 0\;.
\]\end{thm}

\medskip
This results improves on \cite{MR2152502} since it holds as soon as $\lambda_1^\infty=0$. In the next section, we shall use a similar method to interpolate $\Ip$ in terms of $\Ep$ and $\Kp$.

\section{Generalized entropies and nonlinear diffusion equations}

\subsection{Generalized entropies}\label{Sec:GeneralizedEntropies}

Assume now that $v$ is a nonnegative solution of the {\sl weighted porous media equation\/}
\be{Eqn:PM}
v_t=\Dg v^m\;.
\ee
We refer for instance to \cite{0661} or \cite{Bakry02} for such an equation. Again assume that $d\gamma$ is a probability measure and define for any $p\in (1,2)$ the entropy
\[
\frac 1{m+p-2}\ig{\Big[v^{m+p-1}-M^{m+p-1}-(m+p-1)M^{m+p-2}\,(v-M)\Big]}
\]
where $M=\ig {v(x,t)}$ does not depend on $t>0$. The results in \cite{BBDGV2,BBDGV-Cras} suggest that the method could also work even if $d\gamma$ is not bounded. Taking into account the homogenity and up to a change of the time scale, we may assume that $M=1$ and consider
\begin{eqnarray*}
\Epm(t)&:=&\frac 1{m+p-2}\ig{\Big[v^{m+p-1}-1\Big]}\;,\\
\Ipm(t)&:=&c(m,p)\ig{|\D s|^2}\;,\\
\Kpm(t)&:=&\ig{s^{\beta(m-1)}\,|\Dg s|^2}+\alpha\!\ig{s^{\beta(m-1)}\,\Dg s\,\frac{|\D s|^2}s}\;,
\end{eqnarray*}
with $v=:s^\beta$, $\beta:=(p/2+m-1)^{-1}$, $\alpha:=\beta\,m-1=\frac{2-p}{p+2(m-1)}$ and
\[
c(m,p)=\frac{4\,m\,(m+p-1)}{(2m+p-2)^2}\;.
\]
Written in terms of $s$, the evolution is governed by
\be{Eqn:PM2}
\frac 1m\,s_t=s^{\beta(m-1)}\left[\Dg s+\alpha\,\frac{|\D s|^2}s\right]
\ee
since
\[
\Dg s^q=q\,s^{q-1}\left(\Dg s+(q-1)\,\frac{|\D s|^2}s\right)\,.
\]
Using \eqn{Eqn:PM} and \eqn{Eqn:PM2}, a computation shows that
\begin{eqnarray*}
\frac d{dt}\Epm(t)&:=&-\Ipm(t)\;,\\
\frac 1m\,\frac d{dt}\Ipm(t)&:=&-2\,c(m,p)\,\Kpm(t)\;.
\end{eqnarray*}

\subsection{A functional inequality}

The next step of the method is to establish a functional inequality which replaces~\eqn{Eqn:GPI}. In case $m=1$, we had $\Kp\geq p\,\lambda_1(p)\,\Ip/4$ if $\lambda_1(p)>0$ (see Sections~\ref{Sec:NonLoc}). Here we find a more complicated expression which involves $\Epm$, $\Ipm$ and~$\Kpm$, in the spirit of the computations of Section~\ref{Sec:Refined}.

Using
\be{Eqn:Pwrq}
s^{q-1}\,\Dg s=\frac 1q\,\Dg s^q-(q-1)\,s^{q-2}\,|\D s|^2=\frac 1q\,\Dg s^q-\frac 4{q^2}\,(q-1)\,|\D s^{q/2}|^2
\ee
with $2(q-1)=\beta(m-1)$, we can write
\begin{eqnarray*}
q^2\ig{s^{\beta(m-1)}\,|\Dg s|^2}=\ig{|\Dg s^q|^2}& +&\left(\frac 4q\,(q-1)\right)^2\ig{|\D s^{q/2}|^4}\\& -& 2\cdot\frac 4q\,(q-1)\ig{\Dg s^q\,|\D s^{q/2}|^2}\;.
\end{eqnarray*}
Assume that $\Omega$ is convex. By \eqn{Ineq:IPP1}, we get
\[
\ig{|\Dg s^q|^2} \geq \ig{|\DD s^q|^2}+\ig{\DD F\,\D s^q\cdot\D s^q}\;,
\]
\begin{eqnarray*}\hspace*{-24pt}
q^2\ig{s^{\beta(m-1)}\,|\Dg s|^2}&\geq&\ig{|\DD s^q|^2}+\ig{\DD F\,\D s^q\cdot\D s^q}\\
&&+\;\left(\frac 4q\,(q-1)\right)^2\ig{|\D s^{q/2}|^4}\\
&&-\;2\cdot\frac 4q\,(q-1)\ig{\Dg s^q\,|\D s^{q/2}|^2}\;.
\end{eqnarray*}
With the notation $z:=s^{q/2}$, this can be rewritten as
\begin{eqnarray}\hspace*{-24pt}
q^2\ig{s^{\beta(m-1)}\,|\Dg s|^2}&\geq&\ig{|\DD s^q|^2}+\ig{\DD F\,\D s^q\cdot\D s^q}\nonumber\\
&&+\;\left(\frac 4q\,(q-1)\right)^2\ig{|\D z|^4}\label{Eqn:PMa}\\
&&-\;2\cdot\frac 4q\,(q-1)\ig{\Dg s^q\,|\D z|^2}\;.\nonumber
\end{eqnarray}

\medskip On the other hand, by \eqn{Eqn:Pwrq},
\begin{eqnarray}
&&\hspace*{-24pt}\ig{s^{\beta(m-1)}\,\Dg s\,\frac{|\D s|^2}s}\nonumber\\
&=&\ig{s^{q-1}\,\Dg s\cdot s^{q-2}\,|\D s|^2}\nonumber\\
&=&\ig{\left(\frac 1q\,\Dg s^q-\frac 4{q^2}\,(q-1)\,|\D s^{q/2}|^2\right)\cdot \frac 4{q^2}\,|\D s^{q/2}|^2}\nonumber\\
&=&\frac 4{q^2}\ig{\left(\frac{\Dg s^q}q\,|\D z|^2-\frac 4{q^2}\,(q-1)\,|\D z|^4\right)}\label{Eqn:PMb}
\end{eqnarray}
where we again used the identity $2(q-1)=\beta(m-1)$.

\medskip Collecting the estimates \eqn{Eqn:PMa} and \eqn{Eqn:PMb}, we get
\begin{eqnarray*}\hspace*{-24pt}
\Kpm&\geq&\frac 1{q^2}\left[\;\ig{|\DD s^q|^2}+\ig{\DD F\,\D s^q\cdot\D s^q}\right]\\
&&\qquad+\;\frac 1{q^2}\,\left({\textstyle \frac 4q\,(q-1)}\right)^2\ig{|\D z|^4}\\
&&\qquad-\;\frac 2{q^2}\,\cdot\frac 4q\,(q-1)\ig{\Dg s^q\,|\D z|^2}\\
&&\qquad+\;\frac 4{q^2}\,\alpha\,\ig{\left(\frac{\Dg s^q}q\,|\D z|^2-\frac 4{q^2}\,(q-1)\,|\D z|^4\right)}\;,
\end{eqnarray*}
\begin{eqnarray}\hspace*{-24pt}
\Kpm&\geq&\frac 1{q^2}\left[\;\ig{|\DD s^q|^2}+\ig{\DD F\,\D s^q\cdot\D s^q}\right]\nonumber\\
&&\qquad+\;16\,\frac{q-1}{q^4}\,(q-1-\alpha)\ig{|\D z|^4}\nonumber\\
&&\qquad+\;4\,\frac{\alpha-2\,q+2}{q^3}\ig{\Dg s^q\,|\D z|^2}\;.\label{Eqn:PMc}
\end{eqnarray}

\medskip By definition of $z:=s^{q/2}$, we get
\[
\Dg (s^q)=2\,|\D z|^2+2\,s^{q/2}\,\Dg (s^{q/2})\;,
\]
so that
\[
\ig{\Dg (s^q)\,|\D z|^2}=2\ig{|D z|^4}+2\ig{s^{q/2}\,\Dg (s^{q/2})\,|\D z|^2}\;.
\]
Using integration by parts,
\begin{eqnarray*}
&&\hspace*{-24pt}\ig{s^{q/2}\,\Dg (s^{q/2})\,|\D z|^2}\\&=&-\ig{|D z|^4}-2\ig{s^{q/2}\,\DD (s^{q/2})\,:\,\D z\otimes\D z}\\
&=&-\ig{|D z|^4}-\ig{\left[\DD(s^q)-2\,\D z\otimes\D z\right]\,:\,\D z\otimes\D z}\\
&=&\ig{|D z|^4}-\ig{\DD(s^q)\,:\,\D z\otimes\D z}
\end{eqnarray*}
since
\[
\frac 12\,\DD(s^q)=s^{q/2}\,\DD (s^{q/2}) + \D z\otimes\D z\;.
\]
Hence
\[
\ig{\Dg (s^q)\,|\D z|^2}=4\ig{|D z|^4}-2\ig{\DD(s^q)\,:\,\D z\otimes\D z}\;.
\]
Summarizing, with $\mathsf b=\frac 8{q^3}\,(\alpha+2-2\,q)$ and $\mathsf c=\frac{16}{q^4}\,(q-1)\,(q-1-\alpha)+2\,b$, \eqn{Eqn:PMc} becomes
\begin{eqnarray}\hspace*{-24pt}
\Kpm&\geq&\frac 1{q^2}\left[\;\ig{|\DD s^q|^2}+\ig{\DD F\,\D s^q\cdot\D s^q}\right]\nonumber\\
&&\qquad+\;\mathsf c\ig{|\D z|^4}-\;\mathsf b\ig{\DD(s^q)\,:\,\D z\otimes\D z}
\;.\label{Eqn:PMc2}
\end{eqnarray}

\medskip Let
\[
V(x):=\inf_{\xi\in S^{d-1}}\big(\DD F(x)\,\xi,\xi\big)
\]
and define for any $\theta\in (0,1)$
\[
\lambda_1(m,\theta):=\inf_{w\in H^1(\Omega,d\gamma)\setminus\{0\}}\frac{\ig{\Big((1-\theta)\,|\D w|^2+V\,|w|^2\,\Big)}}{\ig{|w|^2}}\;.
\]
Assume that for some $\theta\in (0,1)$,
\be{Hyp:MP}
\lambda_1(m,\theta)>0\;.
\ee
{}From \eqn{Eqn:PMc2}, we get
\begin{multline*}
\Kpm\geq\,\frac{\lambda_1(m,\theta)}{q^2}\ig{|\D s^q|^2}+\;\mathsf a\ig{|\DD s^q|^2}-\;\mathsf b\ig{\DD(s^q)\,:\,\D z\otimes\D z}\\
+\;\mathsf c\ig{|\D z|^4}\;.
\end{multline*}
where $\mathsf a=\mathsf a(\theta)=\frac\theta{q^2}$. A tedious but elementary computation shows that the condition $\mathsf b^2-4\,\mathsf a(\theta)\,\mathsf c<0$ amounts to
\[
(p+2\,m-4)^2+\Big[5\,m^2+2\,(2\,p-7)\,m+(p-3)^2\Big]\,\theta<0\;,
\]
which defines in terms of $(m,p)$ a bounded set $\mathsf E_\theta$ such that $\partial\mathsf E_\theta$ is an ellipse. The set $\mathsf E_1$ is centered at $(1,3/2)$, contained in the rectangle $\big[1-\sqrt 2/2,1+\sqrt 2/2\big]\times\big[0,3\big]$ and
\[
\mathsf E_1\cap \{(m,p)\,:\,m=1\}=\{1\}\times(1,2)\;.
\]
See Figure 2.1. As a function of the parameter $\theta\in (0,1)$, the family $\mathsf E_\theta$ is decreasing and $\lim_{\theta\to 0_+}\bar{\mathsf E}_\theta=\{(m=1,\,p=2)\}$.

Coming back to the estimate of $\Kpm$, if $(m,p)\in\mathsf E_\theta$, then
\[
\Kpm\geq\kappa_1\ig{|\D s^q|^2}+\kappa_2\ig{|\D z|^4}
\]
with
\[
\kappa_1:=\frac{\lambda_1(m,\theta)}{q^2}\quad\mbox{and}\quad\kappa_2:=c-\frac{\mathsf b^2}{4\,\mathsf a(\theta)}\;.
\]
The next step is based on the interpolation of $\ig{|\D s|^2\!}$ between $\ig{|\D s^q|^2\!}$, $\ig{|\D z|^4\!}=\ig{|\D s^{q/2}|^4\!}$ and $(m+p-2)\,\Epm+1=\ig{s^{1+\frac p{p+2(m-1)}}\!}$. Written in terms of $q$, we observe that $(m+p-2)\,\Epm+1=\ig{s^{2(2-q)}\!}$.
By applying Young's inequality with exponents $4$, $2$ and $4$, we get
\begin{multline*}
\frac{q^2}2\ig{|\D s|^2}=\ig{\sqrt A\,B^\frac 14\,s^{\frac 12(4-3q)}\;\frac{\D s^q}{\sqrt A}\cdot\frac{\D z}{B^\frac 14}}\\
\leq\frac{A^2\,B}{4}\ig{s^{2(4-3q)}}+\frac{1}{2\,A}\ig{|\D s^q|^2}+\frac{1}{4\,B}\ig{|\D z|^4}
\end{multline*}
for any $A$, $B>0$. Let $\eta>0$ be such that
\[
\kappa_1=\frac{\eta}{q^2\,A}\quad\mbox{and}\quad\kappa_2=\frac{\eta}{2\,q^2\,B}\;.
\]
The inequality becomes
\begin{multline*}
\frac\eta{c(m,p)}\,\Ipm=\eta\,\ig{|\D s|^2}\leq\Kpm+\frac{A^2\,B\,\eta}{4}\ig{s^{2(4-3q)}}\\
\leq\Kpm+\frac{\eta^4}{4\,q^8\,\kappa_1^2\,\kappa_2}\ig{s^{2(4-3q)}}\;.
\end{multline*}
By H\"older's inequality, if $q\in(1,4/3)$, we obtain
\[
\ig{s^{2(4-3q)}}\leq\(\ig{s^{2(2-q)}}\)^\frac{4-3q}{2-q}
\]
since $d\gamma$ is a probability measure. Let
\[
\mathsf K:=\frac 1{4\,q^8\,\kappa_1^2\,\kappa_2}\;.
\]
The condition $q\in(1,4/3)$ can be rewritten as
\[
1<m<p+1\;.
\]
Collecting all these estimates, we end up with the following estimate: for any $\eta>0$,
\be{Ineq:Eta}
\frac\eta{c(m,p)}\,\Ipm\leq\Kpm+\mathsf K\, \Big[(m+p-2)\,\Epm+1\Big]^\frac{4-3q}{2-q}\eta^4\;.
\ee

\begin{figure}[ht]\begin{center}\includegraphics{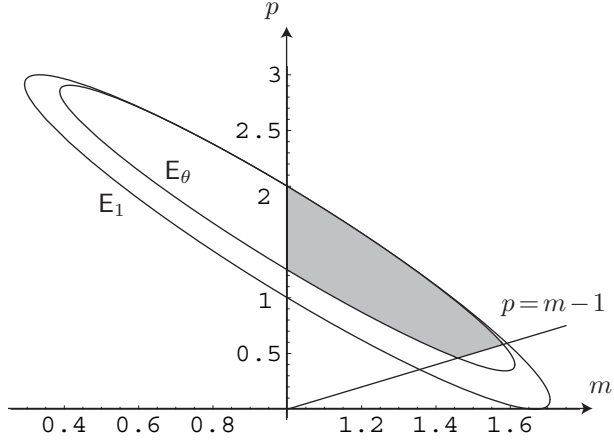}\caption{Admissible parameters $m$ and~$p$ correspond to $(m,p)\in\mathsf E_\theta$, $1<m<p+1$, where the set~$\mathsf E_\theta$ is defined by the condition: $\mathsf b^2-4\,\mathsf a(\theta)\,\mathsf c<0$.}\end{center}\end{figure}

\medskip

\begin{lem}\label{Lem:FunctionalIneq} With the above notations, if $\Omega$ is convex, if \eqn{Hyp:MP} holds for some $\theta\in(0,1)$ and if $(m,p)\in\mathsf E_\theta$ is such that $1<m<p+1$, then
\[\label{eq:functineq}
\Ipm^{\frac{4}{3}}\leq \frac13\,\big[4\,c(m,p)\big]^\frac43\mathsf K^{\frac{1}{3}}\Big[(m+p-2)\,\Epm+1\Big]^{\frac{4-3q}{3(2-q)}}\Kpm\;.
\]\end{lem}

\medskip

\begin{proof} Inequality~\eqn{Ineq:Eta} can be rewritten as
\begin{equation}\label{eq:optim}
f(\eta):=K_1+K_2\,\eta^4-K_3\,\eta\geq 0\quad\forall\;\eta>0\;,
\end{equation}
with
\[
K_1:=\Kpm\;,\quad K_2:=\mathsf K\Big[(m+p-2)\,\Epm+1\Big]^{\frac{4-3q}{2-q}}\;\;\mbox{and}\quad K_3:=\frac{\Ipm}{c(m,p)}\;.
\]
Optimizing on $\eta>0$, we find
\[
0\leq f(\bar\eta)=K_1-\frac{3\,K_3^{4/3}}{4^{4/3}\,K_2^{1/3}}\quad\mbox{where}\quad\bar\eta:=\left(\frac{K_3}{4\,K_2}\right)^{1/3}\,,
\]
which gives the conclusion. \end{proof}

\subsection{Rates of convergence and consequences}

As a consequence of Lemma~\ref{Lem:FunctionalIneq}, we get
\medskip

\begin{thm}\label{Thm:RatePM} If $\Omega$ is convex, if \eqn{Hyp:MP} holds for some $\theta\in(0,1)$ and if $(m,p)\in\mathsf E_\theta$ is such that $1<m<p+1$, then there exists a positive constant $\kappa $, which depends on $\Epm(0)$, such that any smooth solution $u$ of \eqn{Eqn:PM} satisfies, for any $t>0$,
\begin{eqnarray*}
&&\Ipm(t)\leq\frac{\Ipm(0)}{\left[1+\frac\kappa{3}\,\sqrt[3]{\Ipm(0)}\,t\right]^3}\;,\\
&&\Epm(t)\leq\frac{3\,\big[\Ipm(0)\big]^{\frac{8}{3}}}{2\,\kappa\left[1+\frac\kappa{3}\,\sqrt[3]{\Ipm(0)}\,t\right]^2}\;.
\end{eqnarray*}
\end{thm}

\medskip
\begin{proof} If $p>m-1$, then $q\in(1,4/3)$. Using the notations of Section~\ref{Sec:GeneralizedEntropies} and the fact that $\Epm(t)$ is non-increasing, we have
\[
\frac d{dt}\Ipm(t)\leq -\kappa\,\Ipm^{\frac{4}{3}}\;,
\]
with
\[
\kappa:=\frac 32\,m\,\big[4\,\mathsf K\,c(m,p)\big]^{-1/3}\,\Big[(m+p-2)\,\Epm(0)+1\Big]^{-\frac{4-3q}{3(2-q)}}.
\]
An integration from $0$ to $t$ gives the first estimate, and another integration from $t$ to~$\infty$ gives the bound for $\Epm(t)$. \end{proof}

Constants in the above computations are explicit, but not easy to handle. However, we can observe that $\kappa$ can be written as
\[
\kappa=\kappa_0\,\Big[(m+p-2)\,\Epm(0)+1\Big]^{-\frac{4-3q}{3(2-q)}}
\]
where $\kappa_0$ depends only on $m$, $p$, $\theta$ and $\lambda_1(m,\theta)$. Using the estimate for $\Epm$ at $t=0$ in Theorem~\ref{Thm:RatePM}, we obtain the following new functional inequality.
\medskip

\begin{cor}\label{Cor:FunctionalIneq} Under the assumptions of Theorem~\ref{Thm:RatePM}, there exists a positive function $F$ depending only on $m$, $p$, $\theta$ and $\lambda_1(m,\theta)$ such that
\[
F(\Epm)\leq \frac 32\,\Ipm^{\frac{8}{3}}\;,
\]
where $F(s)=O(s)$ as $s\to 0_+$ and $F(s)=O\(s^\frac 2{3(2-q)}\)$ as $s\to\infty$.\end{cor}

\medskip
With the above notations, we can actually write $F$ as
\[
F(s)=\kappa_0\,s\,\Big[(m+p-2)\,s+1\Big]^{-\frac{4-3q}{3(2-q)}}\,,\quad q=\frac{p+3\,(m-1)}{p+2\,(m-1)}\;,
\]
where the dependence of $\kappa_0$ on $m$, $p$, $\theta$ and $\lambda_1(m,\theta)$ can be traced explicitly in the previous computations.

\section{Conclusion}

Our method relies on Conditions \eqn{Hyp:OU} and \eqn{Hyp:MP}. If \eqn{Hyp:OU} is true for some $p=p_0\in[1,2)$, then it is also true for any $p\in[p_0,2]$. Condition \eqn{Hyp:OU} for $p=p_0$ is equivalent to Condition~\eqn{Hyp:MP} for $\theta=\theta_0:=\frac 2{p_0}-1$, and then also holds for any $\theta\in(0,\theta_0)$. The domain of validity of \eqn{Hyp:MP} is represented in Fig. 2.1.

In a linear framework as well as in a non-linear setting corresponding to, respectively, \eqn{Eqn:Heat} and \eqn{Eqn:PM}, the Bakry-Emery method has been adapted to establish functional inequalities which extend the family of generalized Poincar\'e inequalities introduced by W. Beckner. As a consequence, we obtain rates of convergence for a whole family of functionals in each case. It is remarkable that the assumption on the potential can be reduced to the positivity of an eigenvalue, namely Conditions \eqn{Hyp:OU} or \eqn{Hyp:MP}, that is a condition on a linear problem, which turns out to be the same in both cases. Such a condition generalizes the Bakry-Emery criterion, in the sense that the condition is non-local. This improves on the standard criterion. Even in the non-linear case, we recover the full efficiency of the Bakry-Emery method and end up with an inequality relating an entropy functional to its associated Fisher information.

\bigskip\noindent{\it Acknowledgements.} {\small The authors thank the ANR IFO, the Procope project no. 09608ZL and University Paris Dauphine for their support.}

\par\noindent{\scriptsize\copyright~2007 by the authors. This paper may be reproduced, in its entirety, for non-commercial~purposes.}


\bibliographystyle{siam}\small
\bibliography{References.bib}

\medskip\noindent{\small
{\sc J. Dolbeault:}
\address{Ceremade (UMR CNRS no. 7534), Universit\'e Paris-Dauphine, place de Lattre de Tassigny, 75775 Paris C\'edex~16, France.\\
\email{dolbeaul@ceremade.dauphine.fr},
\urladdr{http://www.ceremade.dauphine.fr/$\sim$dolbeaul/}
}\\
{\sc B. Nazaret:}
\address{Ceremade (UMR CNRS no. 7534), Universit\'e Paris-Dauphine, place de Lattre de Tassigny, 75775 Paris C\'edex~16, France.\\
\email{nazaret@ceremade.dauphine.fr},
\urladdr{http://www.ceremade.dauphine.fr/$\sim$nazaret/}
}\\
{\sc G. Savar\'e:}
\address{Universit\`a degli studi di Pavia, Department of Mathematics, Via Ferrata 1, 27100, Pavia, Italy.\\
\email{savare@imati.cnr.it},
\urladdr{http://www.imati.cnr.it/$\sim$savare/}
}}

\medskip\begin{flushright}\today\end{flushright}

\end{document}